\documentclass{amsart}

\advance \topmargin by -\headheight
\advance \topmargin by -\headsep

\usepackage{graphics}
\usepackage{amscd}
\usepackage{amssymb}
\usepackage{longtable}

\newtheorem{thm}{Theorem}[section]
\newtheorem{prop}[thm]{Proposition}
\newtheorem{lem}[thm]{Lemma}
\newtheorem{cor}[thm]{Corollary}

\newtheorem*{thm*}{Theorem}
\theoremstyle{definition}

\newtheorem*{eg}{Example}

\newcommand{\bq}{\mathbb{Q}}
\newcommand{\bz}{\mathbb{Z}}
\newcommand{\br}{\mathbb{R}}
\newcommand{\bc}{\mathbb{C}} 

\newcommand{\coker}{\operatorname{coker}}
\newcommand{\fred}{S^2\times S^2}
\newcommand{\lan}{\left\langle}
\newcommand{\ran}{\right\rangle}
\newcommand{\tet}{\mbox{Tet}}
\newcommand{\oct}{\mbox{Oct}}
\newcommand{\icos}{\mbox{Icos}}
\newcommand{\aut}{\operatorname{Aut}}

\newcommand{\image}{\operatorname{im}}

\newcommand{\fix}{\operatorname{Fix}}
\newcommand{\ind}{\operatorname{Ind}}
\newcommand{\coind}{\operatorname{Coind}}
\newcommand{\cal}{\mathcal}
\newcommand{\gr}{\operatorname{\cal{GR}}}
\newcommand{\coiz}{\coind_{\lan t\ran}^{D_p}(\bz)}
\renewcommand{\hom}{\operatorname{Hom}}
\newcommand{\ad}{\operatorname{Ad}}
\newcommand{\trace}{\operatorname{Trace}}
\newcommand{\matrixb}{\left(\begin{smallmatrix} 0&1\\
	1&0\end{smallmatrix}\right)}

\title{Symmetry groups of four-manifolds}
\author{Michael P. McCooey}
\address{Department of Mathematics and Statistics, McMaster University\\ 
1280 Main Street West\\
Hamilton, Ontario, Canada L8S 4K1}
\date{\today} 
\email{mmccooey@member.ams.org}
\thanks{A portion of this work was supported by a U.S. Department
of Education GAANN fellowship.}

\begin{document}

\begin{abstract}

If a (possibly finite) compact Lie group acts effectively,
locally linearly, and homologically trivially on a closed,
simply-connected four-manifold $M$ with
$b_2(M)\ge 3$, then it  must be isomorphic to a subgroup of
$S^1\times S^1$, and the action must have nonempty fixed-point
set. 

Our results strengthen and complement recent work by
Edmonds, Hambleton and Lee, and Wilczy\'nski, among others.
Our tools include representation theory, finite group theory, and 
Borel equivariant cohomology. 

\end{abstract}

\maketitle

\section{Introduction}

By well-known constructions, any finitely presented group can be
realized as the fundamental group of a closed four-manifold. Any such
group thus acts (by covering translations) on a simply-connected
four-manifold $M$. This action is free, so if $G$ is finite,
the Lefschetz fixed-point theorem implies that it must have a faithful
representation on $H^*(M)$. In contrast, it is natural to ask which groups 
admit homologically trivial, non-free actions on simply-connected 
four-manifolds. The results of this paper provide a nearly complete
answer. Our methods combine cohomological considerations with
local geometric ones,  so in fact our results apply equally
well to manifolds with perfect fundamental groups:
\begin{thm*}[\ref{maintheorem}]
Let $G$ be a (possibly finite) compact Lie group, and suppose
$M$ is a closed four-manifold with $H_1(M; \bz)=0$ and $b_2(M)\ge 2$,
equipped with an effective, locally linear, homologically trivial
$G$-action.
\begin{enumerate}
\item{If $b_2(M)=2$ and $\fix(G)\ne\emptyset$, then $G$ is isomorphic
to a subgroup of $S^1\times S^1$.}
\item{If $b_2(M)\ge 3$, then $G$ is isomorphic to a subgroup
of $S^1\times S^1$, and a fixed point necessarily exists.}
\end{enumerate}
\end{thm*}

In~\cite{MM3}, we treated the case of $\bz_p\times\bz_p$ actions,
and here we apply one result of that paper to study actions of
non-abelian groups which contain a rank two abelian subgroup.

Each simply-connected four-manifold with $b_2(M)\le 1$ admits actions by 
nonabelian groups, as do $\fred$, $\bc P^2\# -\bc P^2$, and 
$\widehat{\bc P^2}\# -\bc P^2$. On the other hand, results of Hambleton and 
Lee show that the only groups which act smoothly and
homologically trivially on
$\bc P^2\#\bc P^2$ are abelian of rank $\le 2$. 
By combining our main theorem with some extra information provided by
the Atiyah-Singer $G$-signature theorem,
we extend this result to the locally linear case.
The question of exactly which nonabelian groups can act 
on each of the other ``small'' manifolds has been treated to some 
extent by various authors, and we summarize their results below.

The main technical tool in our arguments is Borel equivariant cohomology, 
which associates to a  $G$-space $X$ a graded cohomology
module $H_G^*(X)$. If $M$ is a four-manifold and $\Sigma\subset M$
is the set of points with non-trivial isotropy groups, then the
inclusion of $\Sigma$ in $M$ induces an isomorphism $H_G^*(M)\stackrel
{\cong}{\rightarrow} H_G^*(\Sigma)$ in degrees five and higher.
By showing that these modules can not be isomorphic,  we rule
out the possibility of a $G$-action. 
	
To use this observation to prove a theorem of any generality,
we must first make judicious choices of groups $G$ to study, and then
understand the $G$-spaces $M$ and $\Sigma$ sufficiently well to carry
out cohomology calculations. Our choices of $G$ are determined by an
algebraic classification of minimal nonabelian finite groups. The
assumption of homological triviality makes $H_G^*(M)$ relatively
easy to compute. On the other hand, results of Edmonds show that 
$\Sigma$ is a union
of isolated points and  (possibly intersecting) $2$-spheres. Knowing this,
the representation theory of $G$ and its subgroups can be used to 
study the possible arrangements of spheres near any intersection points,
and thus to understand the global structure of $\Sigma$.

Earlier work on the sort of problem we consider includes that of
Edmonds, Hambleton and Lee, and Wilczy\'nski.
Wilczy\'nski~\cite{Wil2} showed in 1987 that the only groups which
can act locally linearly on $\bc P^2$ are the subgroups of $PGL(3, \bc)$,
and his work also applies to manifolds with homology isomorphic to
that of $\bc P^2$. At approximately the same time, Hambleton and
Lee~\cite{HL} proved a similar result for finite groups. Their
1995 paper~\cite{HL2} uses equivariant gauge theory to re-prove this
result for smooth actions and shows much more generally that if $M$ is
a connected sum of $n>1$ copies of $\bc P^2$ and $G$ acts smoothly
and homologically trivially, then $G$ must be abelian of rank $\le 2$.
Edmonds's recent paper~\cite{Actions} inspired the present work,
and includes our main result as a conjecture. 
He shows that if a finite group $G$ acts homologically trivially,
locally linearly, and \emph{pseudofreely} (i.e. with a singular
set consisting only of isolated points) on a simply-connected
four-manifold $M$ with second Betti number at least 
three, then $G$ must be cyclic. We recover his result as a corollary 
of Theorem~\ref{maintheorem}. 
But in the pseudofree case, our proof reduces essentially to his,
and is not a fundamentally different argument.
Note that the  Betti number requirement in that theorem is necessary: 
in~\cite{MM1}, we consider the case $\fred$, where more complicated
groups arise.

\subsection{Structure of the paper:} 

We begin in section~\ref{examples} with some examples and basic 
observations. In the next section, we briefly review the construction
of Borel equivariant cohomology and describe its application to the
problem at hand. Section~\ref{minimalnonabelian} contains a classification
of minimal nonabelian finite groups and provides the framework of our 
argument. They are naturally divided into rank one groups, which have
periodic cohomology, and groups of rank two and higher, whose actions
can be analyzed via their elementary abelian subgroups and the results
of~\cite{MM3}. These two analyses follow, and in the last section we
gather the ingredients and prove the main theorem and two corollaries.

\subsection{Acknowledgments}

An earlier version of this 
work formed part of my Ph.D. dissertation at Indiana University,
and I would like to express my gratitude to my advisor, Allan Edmonds,
for all his help over the years. Two of his papers, \cite{Aspects} and
\cite{Actions}, were especially important influences. Thanks also to Erg\"un
Yal\c{c}in for useful discussions about minimal nonabelian groups.

\section{A few examples and basic results.}\label{examples}

In this section we present a small collection of examples to 
indicate some of the range of possible structures for the singular set 
(denoted in general by $\Sigma$) of
a locally linear group action. We then see that the added
assumption of homological triviality places surprisingly strong
restrictions on  $\Sigma$. We recall two important results
which we shall use often: first, a version of the Lefschetz
Fixed-point Theorem (compare tom Dieck~\cite[page 225]{tomDieck}.):
\begin{thm}\label{Lefschetz}
Let $g: X\rightarrow X$ be a periodic, locally linear map on a
compact manifold $X$. Then 
$$\chi(\fix(g))=\lambda(g) :=\sum_{i=0}^{\dim X}(-1)^i
	\trace{H_*(g)|_{H_i(X, \bq)}}.$$
\end{thm}

Notice in particular that if $g$ acts trivially on homology, 
$\chi(\fix(g))= \chi(X)$. Also recall a theorem of Edmonds~\cite{Aspects}: 
(The original theorem is stated  only for simply-connected manifolds,
but the proof applies just as well whenever $H_1(M;\bz)=0$.)
\begin{thm}\label{H2independent}
Let $\lan g\ran$ be a cyclic group of prime order $p$ which 
acts, preserving orientation, on a closed four-manifold $M$ with $H_1(M)=0$.
If $\fix(g)$ is not purely 2-dimensional,
then the 2-dimensional components of $\fix(g)$ represent independent elements
of $H_2(M; \bz_p)$. If it is purely 2-dimensional, and has $k$ 
$2$-dimensional components, then the $2$-dimensional components span
a subspace of $H_2(M; \bz_p)$ of dimension at least $k-1$, with any
$k-1$ components representing independent elements.
\end{thm}

If the action of $G$ on $M$ is locally linear and preserves orientation,
then the fixed point set of each $g\in G$ will be a submanifold, each
component of which has even codimension. But their dimensions
need not be equal:

\begin{eg}\label{fpdimension} 
Let $n$ be odd. We construct a $\bz_n$ action on $\fred$:
Begin with an action of $g$ on $D^4$
by $g(z, w)=(z, \lambda w)$, where $\lambda = e^{\frac{2\pi i}{n}}$. 
Along the fixed $S^1\subset\partial D^4$, add a 2-handle with the
same action and framing 0, using the equivariant attaching map
$f_0:S^1\times D^2\rightarrow S^3$ which simply sends $(z, w)$ to $(z,w)$.
Next, let $g$ act on another 2-handle via $(z, w)\mapsto (\lambda z, 
\lambda^{-2} w)$. With respect to this action, the framing 2 attaching
map $f_2(z, w)=(z, z^2w)$ is equivariant. Use it to attach the second
2-handle, linking the first once. In the boundary of the quotient, 
attaching this
2-handle amounts to 2-surgery on $S^2\times S^1$, so the result is 
a lens space. Thus the action can be capped off upstairs with a 4-ball
on which $g$ acts linearly, with an isolated fixed point. The resulting
manifold is easily seen to be $\fred$, and the action has a
fixed-point set consisting of two isolated points and a $2$-sphere.

In contrast, Edmonds shows in~\cite{Aspects} that every component of the 
fixed-point set of an involution on a simply connected
spin 4-manifold has the same dimension. Notice also that the fixed 
point set of any cyclic subgroup of $SO(3)\times SO(3)$ is a
product of spheres, so the action we have constructed is not 
equivalent to a ``linear'' one. \end{eg}

\begin{cor} If the singular set of a cyclic group action on $\fred$
consists of two $2$-spheres, the two spheres need not have the same
isotropy groups.
\end{cor}
\begin{proof}
Simply carry out the above construction with a $g$ of order $2n$.
Then $\fix(g)=S^2\cup S^0$, but $\fix(g^n)=S^2\cup S^2$.
\end{proof}

\begin{eg}
Let $\rho$ denote the $\bz_2$ action on $S^2$ of reflection through 
an equator. Then the fixed-point set of the the diagonal $\bz_2$ action 
$(\rho, \rho)$ on $\fred$ is a torus. By taking equivariant connected
sums $\#_{i=1}^g \fred$, we obtain actions whose fixed point sets
are oriented surfaces of any genus.
\end{eg} 

\begin{eg}\label{dihedralexample}Let $\rho$ be as above, 
and let $r_n$ be a rotation of $2\pi /n$ radians
around an axis meeting the equator. Consider actions on $\fred$. 
If $g=(\rho, \rho)$, and $h=(1, r_n)$, then $\fix(g)\cap \fix(h)
\cong S^1\cup S^1$. Also, $\lan g, h\ran\cong D_n$.
\end{eg} 

Homologically trivial actions are simpler: by the  Lefschetz
fixed-point theorem, for each $g\in G$ which acts homologically trivially,
$\chi(\fix(g))=\chi(M)=b_2(M)+2$. By~\cite{Aspects}, $b_1(\fix(g))=0$, 
so the fixed-point set of each element consists of 2-spheres 
and/or isolated points. The intersections of fixed-point sets are
also easy to describe (compare~\cite[4.8]{Assadi}, \cite[2.3]{HL}):
 
\begin{prop}\label{fpintersection} Let $G$ act locally linearly and 
homologically trivially
on a closed, simply connected 4-manifold $M$ with $b_2(M)\ge 1$.
No isotropy representation can reverse orientation on a singular $2$-plane
in a neighborhood of its fixed point.
Thus the intersection of the fixed-point sets of any $g, h\in G$
consists of a (possibly empty) set of 2-spheres and isolated points. 
\end{prop}

\begin{proof}
It follows from Theorem~\ref{H2independent} that
when $b_2(M)\ge 1$, any
$2$-sphere component of the fixed point set of a cyclic group action
represents a nontrivial element of $H_2(M)$. 
A singular $2$-plane near a fixed point is part
of a singular $S^2$, and if the group action were to reverse orientation
locally, it would send $[S^2]$ to $-[S^2]$.

For the second claim, we
need to rule out the possibility that $\fix(g)\cap \fix(h)$ contains
1-dimensional components. Suppose for a contradiction that $S$ is a
circle component of $\fix(g)\cap\fix(h)$. We may assume $G=\lan g,h\ran$.
By local linearity, $G$ acts preserving orientation
on the linking sphere $S^2$ to $S$, so $G$ is polyhedral. Moreover,
the stabilizer of each point on $ S^2$  is cyclic.
If $G$ itself is cyclic, it has a global fixed point on $S^2$, and 
hence $\fix(G, M)$ is 2-dimensional near $S$.  Otherwise,  for any 
non-cyclic subgroup $H$ of $G$, $S$ is also a component of $\fix(H, M)$. 
Since $D_2$ is a subgroup of each of $\tet$, $\oct$, $\icos$, and
$D_n$, for $n$ even, it suffices to rule out dihedral groups.
But in a $D_n$ action on $S^2$, an involution reverses the 
orientation on the order $n$ axis of rotation, which in turn
reverses the orientation on the fixed sphere of the group 
element of order $n$. 
\end{proof}

This proposition places a strong restriction on the allowable linear
representations of isotropy groups, and eventually on the groups themselves,
as we shall see.

\section{The Borel fibration, equivariant cohomology, and the Borel 
spectral sequence}\label{algtop} 

None of the material in this section is new; we include it to 
fix notation and for the convenience of the reader.
For a more thorough discussion, see~\cite{Bredon},
\cite{tomDieck}, or the original source,~\cite{Borel}.

Let $G$ be a compact Lie group, let $X$ be a $G$-space,
and let $\Sigma$ denote the singular 
set of the group action. 
If $EG$ is a contractible, $G$-CW complex on which $G$ acts freely, then
$BG=EG/G$ is a classifying space for $G$-bundles.  
We form the twisted product $X_G=EG\times_G X = (EG\times X)/G$.
The ``Borel fibering'' $X\rightarrow X_G\stackrel{p}{\rightarrow} 
BG$ is naturally induced by
the projection $EG\times X\rightarrow EG$. A similar construction
applies to the singular set $\Sigma$. The \emph{equivariant cohomology}
$H^*_G(X)$ of the $G$-space $X$ is defined as $H^*(X_G)$ (singular
cohomology). Similarly, 
$H^*_G(X, A) := H^*(X_G, A_G)$ whenever $A$ is a closed, $G$-invariant
subspace of $X$. With these definitions, $H^*_G$ is a cohomology
theory with the usual properties. Moreover, $H_G^*(X)$ inherits an
$H^*(G)$-module structure: When $a\in H^*(BG)$ and $x\in H^*_G(X)$,
we define $ax=p^*(a)\cup x$. 
 
Whenever $X$ is a $G$-space, we can apply the Leray-Serre spectral 
sequence to the fibering $X_G\rightarrow BG$. Thus
$$E^{i,j}_2(X)=H^i(BG; H^j(X))\Rightarrow H^{i+j}(X_G)$$
as a sequence of bigraded differential $H^*(BG)$-algebras. The product in 
$E_2$ corresponds to the cup product in $H^*(BG; H^*(X))$ 
(with local coefficients). In particular, if $X$ is connected,
$E^{i,0}$ is naturally identified via $p^*$ with $H^*(BG)$. 
The products and differentials in the spectral sequence 
$E(X)$ are compatible with the
$H^*(G)$-module structure of $H^*(X_G)$. For details, see
Whitehead~\cite[XIII.8]{Whitehead}.
There are appropriate versions of all of this 
with pairs $(X, A)$ whenever $A$ is closed and $G$-invariant.
Henceforth we shall refer to the Leray-Serre spectral sequence of
the Borel fibration simply as the \emph{Borel spectral sequence}.

The following observation
is made in~\cite{Aspects} and~\cite{Actions} in certain forms,
but we include a short proof here for convenience. 

\begin{lem}
Let $M$ be a four-manifold with a locally linear $G$-action, and let 
$\Sigma$ be its singular set. Then
$H^n_G(M)\cong H^n_G(\Sigma)$ for $n>4$.
\end{lem}
\begin{proof}
It follows from  results of 
Spanier~\cite[6.6.2 and 6.9.5]{Spanier} that $\Sigma$ is tautly
embedded in $M$. Thus  we have an isomorphism $\varinjlim
H^*_G(M, U)\cong H^*_G(M, \Sigma)$, as $U$ varies over invariant 
neighborhoods of $\Sigma$. 
Now consider the projection $M\times EG\rightarrow M$.
If we quotient by the $G$-action, we obtain a map $p: M_G\rightarrow M/G$
such that for any $x\in M/G$, $p^{-1}(x)\cong BG_x$. In particular,
when restricted to the complement of the singular set, $p$ 
becomes a fibration with contractible fiber. A trivial application
of the relative spectral sequence of this fibration then shows that
$H^n((M-\Sigma)_G, (U-\Sigma)_G)\cong H^n((M-\Sigma)/G, (U-\Sigma)/G)$ 
for all $n$, whenever $U$  is an invariant open  neighborhood of $\Sigma$. 
By excision again, the latter group is isomorphic to $H^n(M/G, U/G)$.
But since $M$ is four-dimensional, $H^n(M/G, U/G)$ is trivial for 
$n>4$. The lemma follows from the long exact cohomology sequence of the 
pair $(M, A)$ and tautness.
 \end{proof}

To exploit the  isomorphism $H^n_G(M)\cong H^n_G(\Sigma)$,
we need to better understand the topology of the singular set as a $G$-space
and study certain  aspects of the spectral sequences in detail.

Let $j:X\rightarrow EG\times_G X$ be the inclusion of a typical fiber
into $X_G$, and let $A$ be a closed, $G$-invariant subspace of $X$. 
According to tom Dieck~\cite[III.1.18]{tomDieck}, we have:

\begin{prop}
Suppose $H^*(X, A;R)$ is a finitely generated, free $R$-module. 
If $G$ acts trivially on $H^*(X, A)$ and the Borel spectral sequence
$E(X, A)$ collapses, then $H_G^*(X, A)$ is a free $H^*(BG)$-module.
Any set $(x_{\nu}|\nu\in J)$, $x_{\nu}\in H^*_G(X, A)$, such that 
$(j^*x_{\nu}|\nu\in J)$ is an $R$-basis of $H^*(X, A)$, can be taken as 
an $H^*(BG)$-basis. 
\end{prop}

If $M$ has cohomology in only even dimensions and $G$, only in odd
dimensions, then the spectral sequence will collapse automatically.

\begin{cor}\label{freemodule}
Suppose $G$ acts locally linearly
and homologically trivially on a closed  four-manifold
$M$ with $H_1(M)=0$. If $H^*(G)$ vanishes in odd dimensions, 
then $H^*_G(M)$ is a free $H^*(G)$ module on $b_2(M)+2$ generators
corresponding to generators for $H^*(M)$.
\end{cor}

We will also need to consider the spectral sequence for the singular
set $\Sigma$.  It suffices to 
consider, one at a time, the subspaces $G X$, where $X$ is a 
path-component of $\Sigma$. Such a subspace will henceforth be referred
to as a {\em $G$-component} of $\Sigma$.

\section{Minimal nonabelian groups}\label{minimalnonabelian}

In this section we classify the finite non-abelian groups
of which every proper subgroup is abelian. The task of 
classification is not as difficult as one might guess. In 
particular, by requiring that every proper subgroup be abelian, 
rather than only the proper normal subgroups, we bypass 
algebraic questions about simple groups. We begin with a lemma:

\begin{lem}
Let $K$ be a finite abelian $p$-group, and suppose there exists an 
automorphism $\sigma$ of $K$ of prime order $q\ne p$, so that
the resulting representation of $\bz_q$ on any invariant proper
subgroup of $K$ is trivial. Then $K$ has exponent $p$.
\end{lem}

\begin{proof}
The group operation in K will be written additively.
Write $K\cong \bz_{p^{n_1}}\times\dots\times\bz_{p^{n_k}}$, where
$n_1\ge n_2\ge\dots\ge n_k$, and let $a_1, \dots, a_k$ be generators
of the factors in this decomposition.

Let $S=\{x\in K\ |\ \text{order}(x) < p^{n_1}\}$. Observe that $S$ must
be invariant under $\sigma$, so $\sigma |_S$ is trivial. If $n_1=1$,
then $K$ has exponent $p$, and the proof is complete. So assume $n_1>1$.
In this case, $\{x\in K\ |\ px=0\}\subseteq S$.

Write $\sigma(a_1)=\alpha_1 a_1 + \alpha_2 a_2\ + \dots +\alpha_ka_k$,
where $\alpha_i\in \bz_{p^{n_i}}$. Then $\sigma(pa_1)=
p\alpha_1a_1 + \dots + p\alpha_k a_k$.  On the other hand, $pa_1\in S$, 
so $\sigma(pa_1)=pa_1$, and $p\alpha_2a_2=\dots =p\alpha_ka_k=0$,
so we have $\sigma(a_1)=\alpha_1a_1 +x$, where $x\in S$ and $\alpha_1\equiv
1\pmod{p^{n_1-1}}$. 

Now, $\sigma^q(a_1) = \alpha_1^qa_1 + (\alpha_1^{q-1}+ \dots + 
\alpha_1 +1)x$. But $\sigma$ has order $q$, so $\sigma^q(a_1)=a_1$.
Hence $\alpha_1^q\equiv 1\pmod{p^{n_1}}$.

Since $\alpha_1\equiv 1\pmod{p^{n_1-1}}$, we have $\alpha_1 =mp^{n_1-1}+1$
for some $m$. Then by the binomial theorem, 
$$\alpha_1^q = \sum_{i=0}^{q}\binom{q}{i}(mp^{n_1-1})^i\equiv 1
\pmod{p^{n_1}}.$$ 
But $p^{n_1}$ divides all the terms with $i>1$, and the $i=0$ term
is $1$. So $p^{n_1}|qmp^{n_1-1}$, and then $p|m$. Hence $\alpha_1\equiv 1
\pmod{p^{n_1}}$. Now $\sigma^q(a_1)=a_1+qx=a_1$, so $qx=0$. Since
$(p,q)=1$, $x=0$, so $\sigma(a_1)=a_1$. The same argument applies to any 
other element of order $p^{n_1}$, so $\sigma$ is trivial on all of $K$
and hence cannot have order $q$.
   
\end{proof}

\begin{prop}\label{badgroupslist}
Let $G$ be a finite nonabelian group, every proper subgroup of which
is abelian. Then $G$ is one of:
\begin{enumerate}
\item{A minimal nonabelian $p$-group}
\item{$(\bz_p\times\dots\times\bz_p)\rtimes\bz_{q^n}$,}
\end{enumerate}
where $p$ and $q$ are distinct primes, and $n\ge 1$.
\end{prop}

\begin{proof}
Recall a theorem of Burnside (See~\cite[th 7.50]{Rotman}): If $Q$
is a Sylow $q$-subgroup of a finite group $G$ such that 
$Q\subset Z(N_G(Q))$, then
$Q$ has a normal complement $K$. 

Assume $G$ is not a $p$-group. Then every Sylow subgroup is proper,
and hence abelian. If each Sylow subgroup is normal, then $G$ must be 
abelian. So suppose $Q$ is a Sylow $q$-subgroup which is not normal
in $G$. Then $N_G(Q)$ must be abelian, so $Q\subset Z(N_G(Q))$, and
Burnside's theorem gives us a normal complement $K$ such that
$G\cong K\rtimes Q$. Observe that for any other $p\ne q$, a Sylow
$p$-subgroup will be contained in $K$ and hence be normal in $G$.

By minimality, $K$ must be a $p$-group. Also by minimality, $Q$ must be
cyclic, say of order $q^n$, and the semidirect product automorphism
$\sigma$ must have order $q$. One more application of minimality
shows that $\sigma$ must be trivial on any proper submodule of $K$,
and then the lemma applies.
\end{proof}

\begin{eg}
 
It is natural to wonder about rank restrictions on the subgroup $K$
in case 2. Rank one examples exist whenever $q|(p-1)$.
$\tet\cong (\bz_2\times\bz_2)\rtimes\bz_3$ is a familiar rank two
example. More generally, note that $|\aut(\bz_p\times\bz_p)|=(p^2-1)(p^2-p)
=(p-1)^2p(p+1)$, and so an order $q$ automorphism $\sigma$ exists
for any prime $q$ dividing $p-1$ or $p+1$. Moreover, if $q$ is an odd
prime dividing $p+1$, it will not divide $p-1$, and hence a $\sigma$
of order $q$ cannot restrict to a nontrivial automorphism of a 
cyclic subgroup of $\bz_p\times\bz_p$. So the resulting groups
$(\bz_p\times\bz_p)\rtimes\bz_q$ give plenty of examples with rank 
two kernel.

This simple argument can also be used to produce examples with higher rank
kernel.  Notice that $|\aut((\bz_p)^k)| = (p^k-1)(p^k -p)\cdots (p^k-p^{k-1})$,
and the first term of this product factors further as 
$(p-1)(p^{k-1}+ p^{k-2} + \dots +1)$. So whenever $\sum_{i=0}^{k-1}p^i$
has a prime factor which does not divide $p-1$ or any smaller 
$\sum_{i=0}^{k'-1}p^i$ there will be a minimal nonabelian group of type
2 whose kernel $K$ has rank $k$. $(\bz_3\times\bz_3\times\bz_3)\rtimes\bz_{13}$
is an example with $k=3$. In any case, only the ranks 1
and 2 will concern us in applications.
\end{eg}
\smallskip

The minimal nonabelian $p$-groups are classified. According to
Yagita~\cite{Yagita}, who in turn cites Redei~\cite{Redei}, they 
are of two types when $p$ is an odd prime (This classification
can actually be proven using arguments similar to those in the 
proofs above):

\noindent\emph{Type 1.} $G_1(m, n, p)=\lan a, b\ |\ 
 a^{p^m}=b^{p^n}=1,
[a,b]=a^{p^{m - 1}}, m\ge 2\ran$. Thus $G_1\cong \lan a\ran\rtimes
\lan b\ran$, with a semidirect product automorphism $\sigma (a) = b^{-1}ab
=aa^{p^{m -1}}$. 

\noindent\emph{Type 2.} $G_2(m, n, p)=
\lan a,b,c\ |\ a^{p^{m}}=b^{p^{n}}=c^p=1,
[a,b]=c, [c,a]= [c,b] =1\ran$. Observe that when one of $m$ or $n$
(say $m$) is greater than 1, $G_2$ contains a rank three abelian
subgroup generated by $a^p$, $b$, and $c$.  $G_2(1, 1, p)$ is of the
form $\lan a, c\ran\rtimes\lan b\ran$, where $\sigma(a) =ac$. 

When $p=2$, $G_1(2, 1, 2)\cong G_2(1, 1, 2)\cong D_4$, and 
there is exactly one more, the  usual quaternion group 
$D_2^*=\lan a,b\ |\ a^4=1, a^2=b^2, [a,b]=a^2 \ran$.

For future reference, it is useful to organize the minimal nonabelian
groups according to the rank of their elementary abelian subgroups:
\begin{enumerate}
\item{In rank 1, we have the groups $\bz_p\rtimes\bz_{q^n}$ and the
quaternions, $D_2^*$. These groups have periodic cohomology.
Of the rank one groups, the dihedral groups
$D_p$, with $p$ an odd prime, will play a special role for us, as their
irreducible real representations are two-dimensional and not free.}
\item{In rank 2, the groups take the form $(\bz_p\times\bz_p)\rtimes\bz_{q^n}$,
$\bz_{p^m}\rtimes\bz_{p^n}$ (i.e. $G_1(m,n,p)$, when $m\ge 2$),
or $(\bz_p\times\bz_p)\rtimes\bz_p$. The important fact here is that
in each case, there is a \emph{normal} rank two subgroup.}
\item{In rank 3 and higher, the precise structure will be unimportant for us.}
\end{enumerate}

\section{The rank one case}\label{rankone}

The quaternion group $D_2^*$ is a minimal nonabelian group of
rank one; according to our classification, all of the others have
the form $G=\bz_p\rtimes\bz_{q^n}$. Fix generators
$a$ and $b$ for the $\bz_p$ and $\bz_{q^n}$ factors, respectively,
and recall that $\sigma$ denotes the automorphism of $\bz_p$ induced
by conjugation by $b$.  
                               
We will need an explicit calculation of the cohomology groups of $G$.
Recall that the integral cohomology of 
$\bz_p$ is generated by a 
class $t\in H^2(\bz_p; \bz)$. In fact, if $\beta :H^1(\bz_p; \bz_p)\rightarrow 
H^2(\bz_p; \bz)$ is the Bockstein homomorphism, and $s(a)=1$,
then $t=\beta(s)$. 

Let $\sigma(a) = ka$, where $k\in \bz_p^*$ and $k^q\equiv 1\pmod p$.
Then $\sigma_*(t)=kt$, so $\sigma_*(t^i)=k^it^i$. In particular, 
the only powers of $t$ fixed by $\sigma_*$ are those which are 
divisible by $q$. Hence
\begin{equation*}
H^0(\bz_{q^n}; H^j(\bz_p))\cong
\begin{cases}
	\bz &\text{if $j=0$,}\\
	\bz_p & \text{if $j\equiv 0\pmod {2q}$ and $i>0$,} \\
	0 & \text{otherwise.}
\end{cases}
\end{equation*}
With this in mind,  the Hochschild-Serre
spectral sequence shows that:
\begin{equation*}
H^i(\bz_p\rtimes\bz_{q^n}; \bz)\cong
\begin{cases}
	\bz & \text{if $i=0$,}\\
	\bz_{q^n} &\text{if $i$ is even, but $q\not |\ i$,}\\
	\bz_p\oplus\bz_{q^n} &\text {if $i\equiv 0 \pmod{2q}$,}\\
	0 & \text {otherwise.}
\end{cases}
\end{equation*}

A similar calculation applies to $D_2^*$, and yields:

\begin{equation*}
H^i(D_2^*; \bz)\cong
\begin{cases}
	\bz & \text{if $i=0$,}\\
	\bz_2\oplus\bz_2 &\text{if $i\equiv 2\pmod{4}$}\\
	\bz_8 &\text {if $i\equiv 0 \pmod{4}$, and $i>0$,}\\
	0 & \text {otherwise.}
\end{cases}
\end{equation*}

Each of the rank one groups has periodic cohomology, and this
periodicity is reflected in the structure of the equivariant
cohomology modules $H_G^*(M)$ and $H_G^*(\Sigma)$. The details
of the arguments vary as $G$ does, so we consider three cases:
\begin{enumerate}
\item{$G=\bz_p\rtimes\bz_q$, with $q>2$.}
\item{$G=\bz_p\rtimes\bz_{q^n}$, with $n>1$; also $G=D_2^*$.}
\item{The dihedral group $G=D_p$.}
\end{enumerate}

\begin{lem}\label{q=2lemma} 
If $G= \bz_p\rtimes\bz_{q^n}\subset SO(4)$, then $q=2$.\end{lem}
\begin{proof}
The subgroup $G_0=\lan a, b^q\ran$ is cyclic, and by the Brauer theorem
on induced characters, each irreducible complex representation takes
the form $\ind_{G_0}^G(V)$, where $V$ is a complex representation
of $G_0$. In particular, each faithful one has dimension $q$. 
Now, if $G\subset SO(4)$,
then it has a representation on $\bc^4=\br^4\otimes\bc$ which splits
as a sum of irreducibles. Since $G$ is not a nontrivial direct
product, at least one must be faithful, so $q\le 4$. And if $q=3$, $G$ 
would necessarily have a  three-dimensional real representation. But the 
finite subgroups of $O(3)$ are well-known, and $\bz_p\rtimes\bz_{3^n}$ 
is not among them.
\end{proof}

Next we recall some terminology and one
basic result about the cohomology of groups. (See~\cite{Brown}
for details.)
Suppose $G$ is a group, $H\subset G$ a subgroup, and $M$
an $H$-module. We write 
$$\bz G\otimes_{\bz H}M = \ind_H^GM$$
and
$$\hom_{\bz H}(\bz G, M) = \coind_H^GM.$$
When $G$ is finite (or more generally, when $(G:H)$ is finite),
$$\coind_H^GM\cong \ind_H^GM\cong\bigoplus_{g\in G/H}gM.$$ 
The following lemma is elementary, but exceedingly useful:

\begin{lem}[Shapiro's Lemma] If $H\subseteq G$ and $M$ is an $H$-module,
then the composition 
$$H^*(G; \coind_H^GM)\stackrel{r^*}{\rightarrow}
H^*(H; \coind_H^GM)\stackrel{\pi_*}{\rightarrow} H^*(H; M)$$
is an isomorphism,  where
$r:H\hookrightarrow G$ is the inclusion, and $\pi :\coind_H^GM\rightarrow M$
is the canonical projection given by $\pi(\varphi)=\varphi(1)$.
\end{lem}

In particular, the Shapiro isomorphism induces an $H^*(G)$-module
structure on $H^*(H)$. Let $a\in H^*(G)$ and $x\in H^*(H)$.
Tracing through the definitions and using the
naturality properties of the cup product, we find that
$ax=r^*(a)\cup x.$

Now suppose $G=\bz_p\rtimes\bz_q$ acts on $M$, with $q>2$. Lemma~\ref{q=2lemma}
shows that $G$ has no fixed points.
If $\lan a \ran$ were to fix a 2-sphere, $G$ would act on it with a fixed 
point, so $\fix(a)$ must consist only of $b_2+2$ isolated points, permuted
freely by $\bz_q$. Let $m=\frac{b_2+2}{q}$.

$\fix(b)$ will contain, say, $n_1$ isolated points and $n_2$
$2$-spheres. Each of these will form part of a free $\bz_p$-orbit. It
follows that
\begin{equation*}
H^j(\Sigma; \bz)\cong
\begin{cases}
\coind_{\lan a \ran}^G(\bz^m)\oplus\coind_{\lan b\ran}^G(\bz^{n_1+n_2})
	&\text {if $j=0$,}\\
\coind_{\lan b\ran}^G(\bz^{n_2}) & \text{if $j=2$,}
\end{cases}
\end{equation*}
whence, by Shapiro's lemma,
\begin{equation*}
H^i(G; H^j(\Sigma))\cong
\begin{cases}
H^i(\bz_p)^m\oplus H^i(\bz_q)^{n_1+n_2} &\text{if $j=0$,}\\
H^i(\bz_q)^{n_2} &\text{if $j=2$.}
\end{cases}
\end{equation*}

The Borel spectral sequence then shows that in degrees $\ge 4$,
$\gr(H^i_G(\Sigma))$ has period 2 as a graded group. On the other hand,
since $H^*(G)$ has period $2q$, Corollary~\ref{freemodule} implies
that $H^*_G(M)$ must have a period which is divisible by $q$.
So $G=\bz_p\rtimes\bz_q$ cannot act as we supposed.

Next, consider $G=\bz_p\rtimes\bz_{q^n}$, with $n>1$. 
According to Wolf~\cite[5.5.10]{Wolf}, every faithful representation of 
$G$ is free, so if $G$ has a fixed point, it must be isolated (and $q=2$
by Lemma~\ref{q=2lemma}).

On the other hand, if some point $x_0$ is fixed by the cyclic subgroup
$G_0$ generated by $ab^q$, then, since $(ab^q)^p\subset
\lan b^q\ran$, $x_0\in \fix(b^q) \setminus\fix(b)$. This is possible only
if $\fix(b^q)$ is $2$-dimensional around $x_0$. Thus $x_0$ is contained
in a $2$-sphere $S$, and $G/\lan b^q\ran\cong \bz_p\rtimes\bz_q$ must
act effectively upon it. Again, this is only possible if $q=2$.

Thus $G=\bz_p\rtimes\bz_{2^n}$, and the singular set consists of, say,
$n_1$ isolated fixed points and $n_2$ $2$-spheres. Each point and each 
sphere is invariant under $G$, so the action of $G$ is trivial on 
$H^*(\Sigma)$.

Since $G$ has cohomology only in even dimensions, the Borel spectral sequence 
collapses for $H_G^*(\Sigma)$, and shows that 

\begin{equation*}
\gr(H^i_G(\Sigma; \bz))\cong
\begin{cases}
(\bz_{2^n})^{n_1+2n_2}\oplus(\bz_p)^{n_2} & \text{if $i=2\pmod{4}, i>4$,}\\
(\bz_{2^n})^{n_1+2n_2}\oplus(\bz_p)^{n_1+n_2} & \text{if $i=0\pmod{4}, i>4$,}
\end{cases}
\end{equation*}
and that
\begin{equation*}
\gr(H^i_G(M; \bz))\cong
\begin{cases}
	(\bz_{2^n})^{b_2+2}\oplus(\bz_p)^2 & \text{if $i=0\pmod{4}, i>4$,}\\
(\bz_{2^n})^{b_2+2}\oplus(\bz_p)^{b_2} & \text{if $i=2\pmod{4}, i>4$.}
\end{cases}
\end{equation*}
 
Comparison of $p$-torsion shows that $n_1+b_2=2$. In particular, if $b_2\ge 3$,
$\bz_p\rtimes\bz_{2^n}$ can not act as we supposed.

The case of $G=D_2^*$ is similar. The same analysis applies to show that
$\Sigma$ would consist of $n_1$ isolated fixed points and $n_2$ spheres,
and we compute:
$$\gr(H^i(\Sigma))\cong\begin{cases}
(\bz_{8})^{n_2}\oplus(\bz_2\times\bz_2)^{n_1+n_2} 
	&\text{if $i\equiv 2$ mod $4$,}\\
(\bz_2\times\bz_2)^{n_2}\oplus(\bz_{8})^{n_1+n_2}
	&\text{if $i\equiv 0$ mod 4.}
\end{cases}$$

$$\gr(H^i(M_{D_2^*}; \bz))\cong\begin{cases}
(\bz_2)^4\oplus(\bz_{8})^{b_2(M)} &\text{if $i\equiv 2$ mod $4$,}\\
(\bz_{8})^2\oplus(\bz_2)^{2b_2(M)} &\text{if $i\equiv 0$ mod $4$.}
\end{cases}$$

When $k\equiv 2\pmod{4}$, we find that
$4+3b_2(M)=2n_1+5n_2$. And when $k\equiv 0\pmod{4}$, we discover 
that $6+2b_2(M)=3n_1+5n_2$.  Combining these equations, we see
again that $n_1+b_2(M)=2$

Finally, we consider the case $G=D_p = \lan s,t |\ s^p=t^2=1,
t^{-1}st=s^{-1}\ran$, when  $p$ is an odd prime. Since
$D_p$ has $2$-dimensional representations, the situation near
a fixed point can be more complicated than in the earlier cases.

By Proposition~\ref{fpintersection}, if the action of $G$
is homologically trivial, then $\fix(s)$ can not be 
2-dimensional near a global fixed point $x_0$.  
Hence the local representation must take the form
$$s\mapsto \begin{pmatrix}r_p & 0 \\0 & (r_p)^k\end{pmatrix},\quad 
t\mapsto \begin{pmatrix}\rho_1 & 0 \\ 0 & \rho_2\end{pmatrix},$$
where $r_p$ is an order $p$ rotation of a 2-plane, $\rho_1$
and $\rho_2$ are 2-plane reflections, and $k\ne 0 \mod p$.  In
particular, $t$ must fix a 2-plane, and the singular set of the $D_p$
action near $x_0$ is the union of the images of $\fix(t)$ under the
powers of $s$. 

The singular set $\Sigma$ may have several components. 
Some, say $n_1$, of them will contain global fixed points. Suppose $X$
is one such. Observe that
any $x\in \fix(t)\cap\fix(s^kt)$ for some $k$ must in fact be
a fixed point for  all of $D_p$, since these two elements together
generate the group. Since $\fix(t)$ consists of a collection of
isolated points and 2-spheres, 
and since around any fixed point, $\bz_p$ cyclically permutes the
spheres, $X$ must be a union of $p$ 2-spheres intersecting in
$m\ge 1$ fixed points $e^0_1, \ldots, e^0_m$ of $D_p$. In one of
the spheres, choose 1-cells $e^1_1, \ldots, e^1_{m-1}$ connecting,
respectively, $e^0_i$ to $e^0_{i+1}$. View the remainder of the 
sphere as a 2-cell. This determines a CW structure on the sphere
which extends equivariantly to a CW structure on all of $X$. 

The resulting $D_p$-equivariant chain complex for $X$ takes the form:
$$0\rightarrow\ind_{\lan t\ran}^{D_p}(\bz)\stackrel{0}{\rightarrow}
(\ind_{\lan t\ran}^{D_p}(\bz))^{m-1}
\stackrel{\partial_1}{\rightarrow}\bz^m\rightarrow 0.$$
Also note that $\hom_{\bz}(\ind_{\lan t\ran}^{D_p}(\bz), M)\cong
\hom_{\bz[\bz_2]}(\bz[D_p], M)\cong \coind_{\lan t\ran}^{D_p}(M)$. 
It follows that, for any coefficient module $M$,  
\begin{eqnarray*}
H^0(X; M) &=& M\\
H^1(X; M) &=& \bigoplus_{m-1}(\coker(M\rightarrow \coind_{\lan s\ran}^{D_p}(M)))\\
H^2(X; M) &=& \coind_{\lan t\ran}^{D_p}(M).
\end{eqnarray*}

For notational convenience, we will abbreviate 
$\coker(M\rightarrow\coind_{\lan t\ran}^{D_p}(M))$ by $Ck(M)$. 
To calculate $H^*(D_p; Ck(\bz))$, consider the coefficient short exact
sequence
$$1\rightarrow\bz\rightarrow\coiz\rightarrow Ck(\bz)\rightarrow 1,$$
and its associated long exact sequence in cohomology.
By Shapiro's Lemma, $$H^*(D_p; \coind^{D_p}_{\lan t\ran}(\bz))
\cong H^*(\lan t\ran; \bz).$$
And with this identification, the inclusion $\bz\hookrightarrow\coiz$
induces the restriction map $r^*:H^*(D_p; \bz)\rightarrow H^*(\lan t\ran; \bz)$.
Thus we have an exact sequence
$$\cdots\rightarrow H^i(D_p; \bz)\stackrel{r^*}{\rightarrow}H^i(\bz_2; \bz)
\rightarrow H^i(D_p; Ck(\bz))\rightarrow H^{i+1}(D_p; \bz)
\rightarrow\cdots.$$

And it follows that
$$H^n(D_p; Ck(\bz))\cong\begin{cases} 
	\bz_p &\text{if $n\equiv 3$ mod 4,}\\
	0 &\text{otherwise.}
	\end{cases}$$
 
Components of $\Sigma$ which do not contain fixed points are simpler.
Aside from the $n_1$ global fixed points, $\fix(s)$ contains, say, $n_2$
$\lan t\ran$-orbits of 2-spheres (containing two spheres each), 
and $n_3$ $\lan t\ran$-orbits of isolated points (containing
two points in each orbit). And 
the rest of $\fix(t)$ contains, say, $n_4$ 2-spheres and $n_5$
isolated points, each of which forms part of an $s$-orbit. Hence
\begin{eqnarray*}
H^0(\Sigma; \bz) & \cong & \bz^{n_1}\oplus (\coind_{\lan s\ran}^{D_p}(\bz))^{n_2+n_3}
	\oplus (\coind_{\lan t\ran}^{D_p}(\bz))^{n_4+n_5} \\
H^1(\Sigma; \bz) & \cong & \bigoplus_{i=1}^{n_1}(Ck(\bz)^{m_i-1})\\
H^2(\Sigma; \bz) & \cong & (\coind_{\lan t\ran}^{D_p}(\bz))^{n_1+n_4}
	\oplus (\coind_{\lan s\ran}^{D_p}(\bz))^{n_2}.
\end{eqnarray*}

Using Shapiro's lemma, we find:
$$H^n(D_p; H^0(\Sigma))\cong\begin{cases}
	0 &\text{if $n$ is odd,}\\
(\bz_2)^{n_1+n_4+n_5}\oplus(\bz_p)^{n_2+n_3} &\text{if $n\equiv 2$ mod 4,}\\
(\bz_2)^{n_1+n_4+n_5}\oplus(\bz_p)^{n_1+n_2+n_3} 
	&\text{if $n\equiv 0$ mod 4,}
\end{cases}$$

$$H^n(D_p; H^1(\Sigma))\cong\begin{cases}
	\bigoplus_{i=1}^{n_1}\bz_p &\text{if $n\equiv 3$ mod 4,}\\
	0 & \text{otherwise.}
\end{cases}$$

$$H^n(D_p; H^2(\Sigma))\cong\begin{cases}
0 &\text{if $n$ is odd,}\\
(\bz_2)^{n_1+n_4}\oplus(\bz_p)^{n_2} &\text{if $n$ is even.}
\end{cases}$$

\begin{lem}
The spectral sequence for $H^*(\Sigma_{D_p}; \bz)$ collapses.
\end{lem}
\begin{proof}
The sequence has nonzero rows $j=0, 1$ and $2$ only, so it suffices to 
show that the differentials $d_2:E_2^{i, 2}\rightarrow E_2^{i+2, 1}$, 
$d_2:E_2^{i, 1}\rightarrow E_2^{i+2, 0}$, and $d_3:E_3^{i, 2}\rightarrow
E_3^{i+3, 0}$ vanish. But $E_2^{i, 2}$ and $E_2^{i+2, 0}$ vanish in odd 
dimensions, and $E_2^{i, 1}$ vanishes in even dimensions. The
conclusion follows.
\end{proof}

\begin{prop}\label{noDp}{If $D_p$ ($p$ odd, prime) acts locally linearly
and homologically trivially on $M$, then $|\fix(D_p)|+b_2(M)=2$.}
\end{prop}

\begin{proof}
Let $G=D_p$. By Corollary~\ref{freemodule},  
$$\gr(H^8(M_G; \bz))\cong \bz_{2p}\oplus (\bz_2)^{b_2(M)}\oplus\bz_{2p}.$$

And since the sequence for $\Sigma_{D_p}$ also collapses, we have
\begin{eqnarray*}
\gr(H^8(\Sigma_{D_p}; \bz)) &\cong& H^8(D_p; H^0(\Sigma))\oplus
	H^7(D_p; H^1(\Sigma))\oplus H^6(D_p; H^2(\Sigma))\\
	&\cong & (\bz_{2p}^{n_1}\oplus \bz_p^{n_2+n_3}\oplus
	\bz_2^{n_4+n_5})\oplus(\bz_p^{\sum (m_i-1)})
	\oplus(\bz_2^{n_1+n_4}\oplus\bz_p^{n_2}).
\end{eqnarray*}
By comparing the $p$-ranks, we find that $n_1+2n_2+n_3+\sum(m_i-1)=2$.
Since the action has a fixed point, $n_1\ge 1$, so we must have
$n_2=0$ and $n_3\le 1$. 
On the other hand, $\chi(\fix(s)) = n_1+4n_2+2n_3$.
By the Lefschetz Fixed-point theorem, $\chi(\fix (s))
=b_2(M)+2$. This is possible only if 
$n_1=n_3=b_2(M)=1$ or $n_1=2$ and $n_3=b_2(M)=0$ -- in other words, 
when  $|\fix(D_p)+b_2(M)=2$.
When $b_2(M)\ge 2$, as we have assumed, there are no actions with
a fixed point.
\end{proof}

Some of the calculation in this section assume that $b_2(M)\ge 1$, but for
these groups the case $b_2(M)=0$ is easily handled by inspection. 
So if a rank 1 minimal nonabelian group acts, then $|\fix(G)|+b_2(M)=2$.

\section{The rank two and higher cases}\label{ranktwo}

Rank 2 groups have more complicated cohomology and representation theory
than rank 1 groups, so a direct analysis along the lines of the one we have
just finished is rather difficult. (The case $G=D_4$ is especially 
interesting.) But as we have observed, the minimal
nonabelian groups of rank $2$ each contain a \emph{normal} subgroup $G_0$
isomorphic to $\bz_p\times\bz_p$. If we denote by $\Sigma_0$ the singular
set of the embedded $\bz_p\times\bz_p$ action, then all of $G$ must act on
$\Sigma_0$. In~\cite{MM3}, we studied $\bz_p\times\bz_p$ actions in detail,
and in particular, we showed:

\begin{prop}
Suppose $M$ is a closed,
topological four-manifold with $b_2(M)\ge 1$ and $H_1(M) = 0$, equipped with
an effective, homologically trivial, locally linear $\bz_p\times\bz_p$
action. Except for fixed-point free actions which exist
in two cases:
\begin{enumerate}
\item{when $b_2(M)=1$,  $p=3$, and the action is pseudofree, or}
\item{when $M$ has intersection form $\matrixb$ and $p=2$,}
\end{enumerate}
the singular set $\Sigma_0$ consists
of $b_2(M)+2$ spheres equipped with rotation actions, intersecting
pairwise at their poles, and arranged into a single closed loop. 
Each sphere represents a primitive class in $H_2(M; \bz)$, and together
these classes generate $H_2(M)$.
\end{prop}

We shall assume in this section that either $b_2(M)\ge 3$, or
$b_2(M)=2$, but the action has a fixed point. In either case, the
lemma determines the structure of $\Sigma_0$. In particular, we can
observe:

If $b_2(M)\ge 2$, then adjacent spheres $S_i$ and $S_{i+1}$
in $\Sigma_0$ must represent different homology classes.  For
$S_i$ intersects $S_{i-1}$ once, while $S_{i+1}$ does not intersect it at
all. So a homologically trivial action which leaves $\Sigma_0$ 
invariant, and fixes a point, must leave the individual spheres
invariant. And if $b_2(M)\ge 3$, then any 
homologically trivial action on $M$ which
leaves $\Sigma_0$ invariant must in fact leave each sphere in $\Sigma_0$
invariant, and fix their intersection points. 
For the $b_2+2$ spheres represent at least $b_2$ different
homology classes.

It follows that the isotropy representation of $G$, which must be
faithful, splits as a sum of two two-dimensional rotation actions.
This is possible only if $G\subset SO(2)\times SO(2)$. Thus actions of
the minimal nonabelian groups of rank two are ruled out,
as are actions of any groups -- even abelian ones -- of rank $\ge 3$.

\begin{eg} 
As a contrast to the constraints which occur when $b_2(M)\ge 3$,
consider the action of $D_4$ on $\bc P^2$ given in homogeneous 
coordinates by $s:[x,y,z]\mapsto [ix, -iy, z]$, $t:[x,y,z]
\mapsto [y,x,z]$. The fixed-point sets of each element are easy to
calculate. If $G_0 =\lan s^2, t\ran$, we 
find that $\Sigma_0$ consists of three
copies of $S^2$ intersecting at their poles. The action of $D_4$
exchanges two of the spheres and leaves the one fixed by $s^2$
invariant. 

Notice that this action can be obtained from a linear
action on $S^4$ by an equivariant blowup at one of the fixed points.
If we perform the same construction at the other fixed point
(with opposite orientation), the result is an action on $\bc P^2\#
\overline{\bc P}^2$, in which $\Sigma_0$ consists of two spheres fixed by
$s^2$ and left invariant by the whole group, and spheres fixed by
$t$ and $s^2t$, respectively, which are exchanged by the action of $s$.
\end{eg}

\section{Proofs of the main results.}

A simple observation will allow us to generalize our results from finite
groups to compact Lie groups:

\begin{lem}\label{Liegrouplemma}
Every nonabelian compact Lie group contains a nonabelian
finite subgroup.
\end{lem}
\begin{proof}
A nonabelian compact Lie group $G$ has a nontrivial Weyl group 
$W=N_G(T)/Z_G(T)$, where $T$ is a maximal torus.  
$N_G(T)$ is itself a compact Lie group whose identity
component $T$ is (of course) abelian. A structure theorem for
such groups (see~\cite[theorem 6.10]{HofmannMorris} or~\cite{DHLee})
states that there is a finite subgroup $E$ of $G$ such that
$N_G(T)=TE$. In particular, there is some $h\in E$ of finite
order $n$, say, which normalizes $T$ but does not centralize it.

Since $\ad(h)$ acts non-trivially on $L(T)$, there is a closed one-parameter
subgroup $\theta_v$ such that $\ad(h)(v)\ne v$. Choose $k\in\image(\theta_v)$ 
of prime  order $p$ which is not fixed by $\ad(h)$. 
Then $k$ and its conjugates under powers of $h$ generate
a subgroup $\bz_p\times\dots\times\bz_p\subset T$ which is normalized by
$h$. Thus $\lan h, k\ran\cong (\bz_p\times\dots\times\bz_p)\rtimes\bz_n$
is a finite, nonabelian subgroup of $G$.
\end{proof}
 
This brings us to our main result:

\begin{thm}\label{maintheorem}
Let $G$ be a (possibly finite) compact Lie group, and suppose
$M$ is a closed four-manifold with $H_1(M; \bz)=0$ and $b_2(M)\ge 2$,
equipped with an effective, locally linear, homologically trivial
$G$-action.
\begin{enumerate}
\item{If $b_2(M)=2$ and $\fix(G)\ne\emptyset$, then $G$ is isomorphic
to a subgroup of $S^1\times S^1$.}
\item{If $b_2(M)\ge 3$, then $G$ is isomorphic to a 
subgroup of $S^1\times S^1$,
and a fixed point necessarily exists.}
\end{enumerate}
\end{thm}
\begin{proof}
If an action by a nonabelian group exists, then by Lemma~\ref{Liegrouplemma},
an action by a minimal nonabelian finite group $G_0$ exists. 
But we have seen that for every such $G_0$ which acts, $b_2(M)+|\fix(G_0)|=2$.
The existence of fixed points for abelian groups follows easily from 
the Lefschetz Fixed Point theorem and the result of~\cite{MM3} used
in Section~\ref{ranktwo} above.
\end{proof}

As a corollary of this theorem we recover the main result 
of~\cite{Actions}:

\begin{cor}[Edmonds]
If a finite group $G$ acts locally linearly, pseudofreely, and
homologically trivially on a closed, simply connected four-manifold
$X$ with $b_2(X)\ge 3$, then $G$ is cyclic and acts semifreely, and the
fixed point set consists of $b_2(X)+2$ isolated points.
\end{cor}

\begin{proof}
By Theorem~\ref{maintheorem}, $G$ must be abelian of rank
at most two, and have a fixed point. But if $G$ has rank two, it
cannot act freely on the linking sphere to the fixed point, and so 
cannot act pseudofreely.
\end{proof}

Since $\bc P^2\#\bc P^2$ and $\bc P^2\#-\bc P^2$ have the same second
Betti number, and both have diagonalizable intersection forms, our methods
so far do not distinguish between them. Yet the latter admits nonabelian
group actions, while the first admits none -- at least, according to Hambleton
and Lee~\cite{HL2}, no smooth actions. One tool which does see the difference
between the two manifolds is the Atiyah-Singer $G$-signature theorem 
(see~\cite{AtiyahSinger3}
for the smooth case, \cite{Wall} for locally linear actions),
and we apply it here to extend Hambleton and Lee's result to the locally
linear setting:

\begin{cor}
Let $G$ be a compact Lie group with a locally linear, homologically
trivial action on a four-manifold $M$ whose integral cohomology is isomorphic
to that of $\bc P^2\#\bc P^2$. Then $G$ is isomorphic to a subgroup
of $S^1\times S^1$.
\end{cor}
\begin{proof}
As before, we suppose $G$ is minimal nonabelian. Theorem~\ref{maintheorem}
implies that the action must be fixed-point-free. We briefly consider 
the possibilities: 
\begin{enumerate}
\item{$G=\bz_p\rtimes\bz_q$, with $q>2$. The arguments of Section~\ref{rankone}
rule out these actions as before, since $H^*_G(M)$ must have period 2 or 4.
}
\item{$G=\bz_p\rtimes\bz_{q^n}$, with $n>1$, or $G=D_2^*$.
As before, we show easily in the first case that $q$ must equal 2.
Choose generators so that $G=\lan a \ran\rtimes\lan b \ran$. Then
$\Sigma=\fix(b^2)$ consists of two spheres on which $b$ acts by an order
two rotation. Denote the fixed points of this rotation by $x_1, \ldots, x_4$,
and the rotation angles of the $b$ action around $x_i$ by $\varphi_i$ and 
$\psi_i$. According to the $G$-signature theorem, $-\sum_i\cot(\varphi_i/2)
\cot(\psi_i/2)=\sigma(g, M)=2$. But at each point, one of $\varphi_i$ or 
$\psi_i$ equals $\pi$, a contradiction. The same argument applies to $D_2^*$.}

\item{The dihedral group $G=D_p=\lan s\ran\rtimes\lan t\ran$. The 
$G$-signature theorem implies that some component $S\subset\fix(t)$ must
have $[S]\cdot[S]\ne 0$. But $sS\subset\fix(sts^{-1})$, and since the
action is homologically trivial, $[sS]\cdot[S]\ne 0$. 
So $\fix(t)\cap\fix(sts^{-1})\ne\emptyset$. Together, these elements generate
$G$, so $\fix(G)\ne \emptyset$, a contradiction.}
\item{$G$ has rank two. Let $G_0\cong\bz_p\times\bz_p\vartriangleleft G$, and
$\Sigma_0$ be as in Section~\ref{ranktwo}. $\Sigma_0$ contains four
spheres; since they cannot all be left individually invariant by $G$,
two opposing spheres are left invariant and two are interchanged by some
$g\in G$. Thus $g$ fixes four points (two on each invariant sphere), while
$g^2$ fixes two spheres. A contradiction follows just as in case 2 above.} 
\end{enumerate}

\end{proof}

\bibliographystyle{abbrv} 
\bibliography{mybiblio}
\end{document}